# ASYMPTOTIC DISTRIBUTIONS OF THE SIGNAL-TO-INTERFERENCE RATIOS OF LMMSE DETECTION IN MULTIUSER COMMUNICATIONS

BY GUANG-MING PAN,[1] MEI-HUI GUO AND WANG ZHOU[2]

*National Sun Yat Sen University, National Sun Yat Sen University and National University of Singapore*

Let $\mathbf{s}_k = \frac{1}{\sqrt{N}}(v_{1k},\ldots,v_{Nk})^T, k=1,\ldots,K$, where $\{v_{ik}, i,k=1,\ldots\}$ are independent and identically distributed random variables with $Ev_{11} = 0$ and $Ev_{11}^2 = 1$. Let $\mathbf{S}_k = (\mathbf{s}_1,\ldots,\mathbf{s}_{k-1},\mathbf{s}_{k+1},\ldots,\mathbf{s}_K)$, $\mathbf{P}_k = \mathrm{diag}(p_1,\ldots,p_{k-1},p_{k+1},\ldots,p_K)$ and $\beta_k = p_k\mathbf{s}_k^T(\mathbf{S}_k\mathbf{P}_k\mathbf{S}_k^T + \sigma^2\mathbf{I})^{-1}\mathbf{s}_k$, where $p_k \geq 0$ and the $\beta_k$ is referred to as the signal-to-interference ratio (SIR) of user $k$ with linear minimum mean-square error (LMMSE) detection in wireless communications. The joint distribution of the SIRs for a finite number of users and the empirical distribution of all users' SIRs are both investigated in this paper when $K$ and $N$ tend to infinity with the limit of their ratio being positive constant. Moreover, the sum of the SIRs of all users, after subtracting a proper value, is shown to have a Gaussian limit.

**1. Introduction.** Consider a symbol synchronous direct sequence code division multiple access (DS-CDMA) system with $K$ users. The discrete-time model for the received signal $y$ in a symbol interval is

$$\mathbf{y} = \sum_{k=1}^{K} x_k \mathbf{s}_k + \mathbf{w}, \tag{1.1}$$

where the $x_k$ is the symbol transmitted by user $k$, $\mathbf{s}_k \in R^N$ is the signature sequence of user $k$ and $\mathbf{w} \in R^N$ is the noise vector with mean zero and covariance matrix $\sigma^2 I$. We also assume that the symbol vector $\mathbf{x} = (x_1,\ldots,x_K)$

Received April 2006; revised July 2006.
[1]Supported by NSC Grant 94-2816-M-110-005 and NSFC Grants 10471135 and 10571001.
[2]Supported in part by Grants R-155-000-035-112 and R-155-050-055-133/101 at the National University of Singapore.

has a covariance matrix $\mathbf{P}$ where $\mathbf{P} = \text{diag}(p_1, \ldots, p_K)$ with $p_k$ being the received power of user $k$, that is, $Ex_k^2 = p_k$ and that the symbol vector is uncorrelated with the noise (more details can be found in [13]).

The engineering goal is to demodulate the transmitted $x_k$ for each user. Assume that the receiver has already acquired the knowledge of the signature sequences. For user $k$, the linear minimum mean-square error (LMMSE) receiver generates an output in a form $\mathbf{a}_k^T \mathbf{y}$ where $\mathbf{a}_k$ is chosen to minimize the mean-squared error

$$(1.2) \qquad E|x_k - \mathbf{a}_k^T \mathbf{y}|^2.$$

The relevant performance measure is the signal-to-interference ratio (SIR) of the estimate (see [13]), which is defined by

$$(1.3) \qquad \beta_k = p_k \mathbf{s}_k^T (\mathbf{S}_k \mathbf{P}_k \mathbf{S}_k^T + \sigma^2 \mathbf{I})^{-1} \mathbf{s}_k, \qquad k = 1, \ldots, K,$$

where $\mathbf{S}_k$ and $\mathbf{P}_k$ are obtained from $\mathbf{S} = (\mathbf{s}_1, \ldots, \mathbf{s}_K)$ and $\mathbf{P}$ by deleting the $k$th column, respectively.

It is difficult to obtain clear engineering insights from (1.3) since it is dependent on the signature sequences. However, if signature sequences are modeled as being random, one may further proceed with the analysis using random matrix theory when the number of users $K$ and the processing gain $N$ approach infinity, that is, suppose

$$\mathbf{s}_k = \frac{1}{\sqrt{N}} (v_{1k}, \ldots, v_{Nk})^T,$$

$k = 1, \ldots, K$, where $\{v_{ik}, i, k = 1, \ldots\}$ are independent and identically distributed (i.i.d.) random variables. Rigorously speaking, if $v_{ik}$ are random variables, then (1.2) should be viewed as a conditional expectation and at this time it is also necessary to assume that the signature sequences are independent of transmitted symbol and noise.

Indeed, considerable progress has been made in this area. For example, Tse and Hanly in [11] derived the asymptotic SIR under MMSE, a decorrelator receiver and a match filter receiver and fluctuations of SIR have subsequently been considered in [10]. Some related results can be found in [13]. Also see [12], and references therein and see the review paper [1] concerning random matrix theory as well.

However, there are still many open problems in this area. For example, Tse and Zeitouni in [10] asked: What is the empirical distribution of the SIR levels of the users across the system? Is this empirical distribution suitable for characterizing the asymptotic distribution of the SIR for a particular user? Is there any type of "weak asymptotic independence" among users? Also, the asymptotic distribution of the sum of all users' SIRs under MMSE has remained unsolved, which has a close connection with another important performance measure, sum mutual information or spectral efficiency (suitable scaling) (for more information concerning the sum mutual information or the spectral efficiency, see [9] and [14]).

In this paper we will answer the above questions. In other words, we will derive the joint asymptotic distribution of the SIRs for different users and the limiting empirical distribution of the SIRs of the users across the system. The sum of the SIRs for all users, after subtracting a proper value, is also shown to have a Gaussian limit, which gives the asymptotic distribution for sum mutual information under MMSE.

Before stating our main results, we will introduce some notation. Write $\mathbf{B}_N = \mathbf{SPS}^T$, whose empirical spectral distribution (ESD) is denoted by $F^{\mathbf{B}_N}$. The ESD of power matrix $\mathbf{P}$ is denoted by $H_N$. Let $c_N = K/N$. $F^{c,H}(x)$ and $H$ will denote the weak limits of the distribution functions $F^{\mathbf{B}_N}$, $H_N$ respectively, as $N, K \to \infty$ if the limits exist. Define $b = \int (x + \sigma^2)^{-1} dF^{c,H}(x)$ and $b_N = \int (x + \sigma^2)^{-1} dF^{c_N, H_N}(x)$, where $F^{c_N, H_N}(x) = F^{c,H}(x)|_{c=c_N, H=H_N}$.

THEOREM 1.1. *Suppose that:*

(a) $\{v_{ij}, i, j = 1, \ldots\}$ *are i.i.d. with* $Ev_{11} = 0$, $Ev_{11}^2 = 1$ *and* $Ev_{11}^6 < \infty$.
(b) $H_N$ *converges weakly to some distribution function* $H$ *and the elements of* $P$ *are bounded by some constant.*
(c) $K/N \to c > 0$ *as* $N \to \infty$.

*Then, for any finite integer* $m$

(1.4) $$(\sqrt{N}(\beta_1 - p_1 b_N), \ldots, \sqrt{N}(\beta_m - p_m b_N)) \xrightarrow{D} N(0, \mathbf{C})$$

*with covariance matrix*

(1.5) $$\mathbf{C} = \left( 2 \int \frac{dF^{c,H}(x)}{(x+\sigma^2)^2} + \left( \int \frac{dF^{c,H}(x)}{x+\sigma^2} \right)^2 (Ev_{11}^4 - 3) \right) \\ \times \operatorname{diag}(p_1^2, \ldots, p_m^2).$$

REMARK 1.1. Theorem 1.1 indicates that the asymptotic independence of the SIR among users holds, as conjectured by Tse and Zeitouni in [10]. This theorem also includes Theorem 4.5 of [10] as a special result. Actually Tse and Zeitouni in [10] only derived the asymptotic distribution for a single SIR under the conditions that $p_1 = \cdots = p_K$ and $v_{11}$ is symmetric.

THEOREM 1.2. *In addition to the assumptions* (b) *and* (c) *of Theorem 1.1, we suppose:* $(a')$ $\{v_{ij}, i, j = 1, \ldots\}$ *are i.i.d. with* $Ev_{11} = 0$ *and* $Ev_{11}^4 < \infty$. *Then*

(1.6) $$G_N(x) \xrightarrow{i.p.} G(x),$$

where

$$G_N(x) = \frac{1}{K} \sum_{k=1}^{K} I(\beta_k \leq x); \tag{1.7}$$

*i.p.* denotes the convergence in probability. Moreover, the Stieltjes transform of $G(x)$ is $\int (bx - z)^{-1} dH(x)$.

REMARK 1.2. Theorem 1.2 characterizes the empirical distribution function of the SIRs for different users, and, simultaneously, it reveals that the asymptotic empirical distribution of the SIRs for a whole system is different from the asymptotic distribution of the SIR for a particular user, which is normally distributed, as shown in Theorem 1.1. For example, consider a special case $p_1 = \cdots = p_K = p$; then one can easily obtain $G(x) = I(pb \leq x < \infty)$.

REMARK 1.3. Indeed, the convergence mode in Theorem 1.2 can be strengthened to converge with probability 1 according to Theorem 7.1 in [4]. In that paper, a more flexible model is employed; they show that the corresponding SIR converges with probability 1 and also provide uniform convergence of the SIRs for all users. It is interesting to consider how to derive the asymptotic distribution of the SIRs under their model.

THEOREM 1.3. *In addition to assumptions* (b) *and* (c) *of Theorem 1.1, suppose that* (a'') $\{v_{ij}, i, j = 1, \ldots\}$ *are i.i.d. with* $Ev_{11} = 0, Ev_{11}^2 = 1$ *and* $Ev_{11}^4 = 3$. (d) $\int x(1 + xb)^{-2} dH(x) = \int x\, dH(x) \int (1 + xb)^{-2} dH(x)$ *and* $\int x^2\, dH(x)(\int (1 + xb)^{-2} dH(x))^2 + \int x^2(1 + xb)^{-4} dH(x) = 2\int x^2(1 + xb)^{-2} dH(x) \int (1 + xb)^{-2} dH(x)$. *Then we have*

$$\sum_{k=1}^{K} (\beta_k - b_N p_k) \xrightarrow{D} N(\mu, \rho), \tag{1.8}$$

where

$$\begin{aligned}
\mu \int \frac{dH(x)}{(1+xb)^2} \\
= 2c \int \frac{dF^{c,H}(x)}{(x+\sigma^2)^2} \int \frac{x^2}{(1+xb)^3} dH(x) \\
- \frac{1}{2\pi} \int \frac{1}{(x+\sigma^2)^2} \arg\left(1 - c\int \frac{t^2 m^2(x)}{(1+tm(x))^2} dH(t)\right) dx
\end{aligned} \tag{1.9}$$

and

$$\begin{aligned}
\rho\left(\int \frac{dH(x)}{(1+xb)^2}\right)^2 \\
= -\frac{1}{2\pi^2} \iint \frac{(d/dz_1)m(z_1)(d/dz_2)m(z_2)}{(z_1+\sigma^2)(z_2+\sigma^2)(m(z_1)-m(z_2))^2} dz_1\, dz_2,
\end{aligned} \tag{1.10}$$

where the contours for $z_1$ and $z_2$ are nonoverlapping and closed and are taken in the positive direction in the complex plane, both enclosing the support of $F^{c,H}(x)$. Here $m(z)$ represents the Stieltjes transform of $F^{c,H}(x)$ and $\Im m(x) = \lim_{z \to x} \Im m(z)$.

REMARK 1.4. Assumption (d) is satisfied when $p_1 = \cdots = p_K = p$ and in this case the formulas (1.9) and (1.10) can be simplified as

$$\frac{\mu}{(1+p)^2} = 2c \int \frac{dF^c(x)}{(x+\sigma^2/p)^2} \int \frac{x^2}{(1+xb)^3} dH(x)$$

(1.11)
$$- \frac{1}{4(a(c)+\sigma^2/p)} - \frac{1}{4(b(c)+\sigma^2/p)}$$

$$+ \frac{1}{2\pi} \int_{a(c)}^{b(c)} \frac{dx}{(x+\sigma^2/p)\sqrt{4c - (x-1-c)^2}}$$

and

(1.12)
$$\rho = (1+p)^4 \frac{2c}{((\sigma^2/p + c - 1)^2 + 4\sigma^2/p)^2},$$

where $a(c) = (1-\sqrt{c})^2$ and $b(c) = (1+\sqrt{c})^2$ and the expression of $F^c$ is referred to [6].

From Theorem 1.3, we can obtain the following corollary concerning sum mutual information under MMSE.

COROLLARY 1.1. *Under the assumptions of Theorem 1.3,*

(1.13)
$$\sum_{k=1}^{K} (\log(1+\beta_k) - \log(1+b_N p_k)) \xrightarrow{D} N(\mu_1, \rho_1),$$

*where*

$$\mu_1 = \mu \int \frac{dH(x)}{1+xb} - c \int \frac{dF^{c,H}(x)}{(x+\sigma^2)^2} \int \frac{x^2}{(1+xb)^2} dH(x),$$

$$\rho_1 = \rho \left( \int \frac{dH(x)}{1+xb} \right)^2.$$

As is seen from Corollary 1.1, the sum mutual information normalized by $N$ (which is the spectral efficiency which is relevant in wireless communications) converges (see also [14]). So one can guess that the small fluctuation, when expanded by a factor of $N$, appears to be Gaussian. However, it is not an easy task to prove it.

The organization of the paper is as follows. Section 2 establishes Theorem 1.1. The proof of Theorem 1.2 is provided in Section 3. The proof of Theorem 1.3 and Remark 1.4 is included in Section 4 and the proof of Corollary 1.1 is contained in the last section. Throughout this paper, $M$ may denote different constants on different occasions and $\|\cdot\|$ denotes the spectral norm of a matrix or the Euclidean norm of a vector. Also, set $\mathbf{A} = \mathbf{SPS}^T + \sigma^2 \mathbf{I}$, $\mathbf{A}_k = \mathbf{S}_k \mathbf{P}_k \mathbf{S}_k^T + \sigma^2 \mathbf{I}$, $k=1,\ldots,K$, to simplify notation.

**2. Proof of Theorem 1.1.** Before beginning with the proof, we first state a lemma.

LEMMA 2.1. *Let $a_{jj}$ be the jth diagonal element of $\mathbf{A}_1^{-1}$. Under the assumptions of Theorem 1.1,*

$$(2.1) \qquad \lim_{N\to\infty} \frac{1}{N} \sum_{j=1}^{N} a_{jj}^2 \xrightarrow{i.p.} b^2.$$

PROOF. From the well-known matrix inverse formula, we have

$$(2.2) \qquad a_{11} = \frac{1}{(\mathbf{S}_1 \mathbf{P}_1 \mathbf{S}_1^T + \sigma^2 \mathbf{I})_{11} - \hat{\mathbf{s}}_1^T \hat{\mathbf{S}}_1^T (\hat{\mathbf{S}}_1 \hat{\mathbf{S}}_1^T + \sigma^2 \mathbf{I})^{-1} \hat{\mathbf{S}}_1 \hat{\mathbf{s}}_1},$$

where

$$\hat{\mathbf{s}}_j = \frac{1}{\sqrt{N}} (\sqrt{p_2} v_{j2}, \ldots, \sqrt{p_K} v_{jK})^T, \qquad \hat{\mathbf{S}}_1^T = (\hat{\mathbf{s}}_2, \ldots, \hat{\mathbf{s}}_N), \qquad j = 1, \ldots, N,$$

and $(\mathbf{S}_1 \mathbf{P}_1 \mathbf{S}_1^T + \sigma^2 \mathbf{I})_{11}$ is defined in (2.3).

Applying the Helly–Bray theorem one can find

$$(2.3) \quad (\mathbf{S}_1 \mathbf{P}_1 \mathbf{S}_1^T + \sigma^2 \mathbf{I})_{11} = \frac{1}{N} \sum_{k=2}^{K} p_k v_{1k}^2 + \sigma^2 \xrightarrow{i.p.} c \int x \, dH(x) + \sigma^2,$$

where we also use the fact that

$$E \left| \frac{1}{N} \sum_{k=2}^{K} p_k (v_{1k}^2 - 1) \right|^2 = \frac{(E v_{11}^2 - 1)^2}{N^2} \sum_{k=2}^{K} p_k \to 0,$$

as $N \to \infty$.

It is observed that

$$\hat{\mathbf{s}}_j = \mathrm{diag}(\sqrt{p_2}, \ldots, \sqrt{p_K}) \tilde{\mathbf{s}}_j, \qquad \hat{\mathbf{S}}_1^T = \mathrm{diag}(\sqrt{p_2}, \ldots, \sqrt{p_K}) \tilde{\mathbf{S}}_1^T$$

with $\tilde{\mathbf{s}}_j = \frac{1}{\sqrt{N}} (v_{j2}, \ldots, v_{jK})^T$ and $\tilde{\mathbf{S}}_1^T = (\tilde{\mathbf{s}}_2, \ldots, \tilde{\mathbf{s}}_N)$. This, together with Lemma 2.7 in [2], implies

$$(2.4) \quad \begin{aligned} E \bigg| \hat{\mathbf{s}}_1^T \hat{\mathbf{S}}_1^T (\hat{\mathbf{S}}_1 \hat{\mathbf{S}}_1^T + \sigma^2 \mathbf{I})^{-1} \hat{\mathbf{S}}_1 \hat{\mathbf{s}}_1 - \frac{1}{N} \mathrm{tr} \mathbf{P}_1^2 \tilde{\mathbf{S}}_1^T (\tilde{\mathbf{S}}_1 \mathbf{P}_1 \tilde{\mathbf{S}}_1^T + \sigma^2 \mathbf{I})^{-1} \tilde{\mathbf{S}}_1 \bigg|^2 \\ \leq \frac{M E v_{11}^4}{N} E \| \mathbf{P}_1 \tilde{\mathbf{S}}_1^T (\tilde{\mathbf{S}}_1 \mathbf{P}_1 \tilde{\mathbf{S}}_1^T + \sigma^2 \mathbf{I})^{-1} \tilde{\mathbf{S}}_1 \mathbf{P}_1 \| \to 0, \end{aligned}$$

as $N \to \infty$, where we use
$$\|\mathbf{P}_1\tilde{\mathbf{S}}_1^T(\tilde{\mathbf{S}}_1\mathbf{P}_1\tilde{\mathbf{S}}_1^T + \sigma^2\mathbf{I})^{-1}\tilde{\mathbf{S}}_1 P_1\|$$
$$\leq M\|\tilde{\mathbf{S}}_1^T(\tilde{\mathbf{S}}_1\mathbf{P}_1\tilde{\mathbf{S}}_1^T + \sigma^2\mathbf{I})^{-1}\tilde{\mathbf{S}}_1\mathbf{P}_1\| \leq M.$$

Let $\breve{\mathbf{s}}_j = (v_{2j}, \ldots, v_{Nj})^T$, $\mathbf{B}_1^{-1} = (\tilde{\mathbf{S}}_1\mathbf{P}_1\tilde{\mathbf{S}}_1^T + \sigma^2\mathbf{I})^{-1}$, $\mathbf{B}_{1j}^{-1} = (\tilde{\mathbf{S}}_{1j}\mathbf{P}_{1j}\tilde{\mathbf{S}}_{1j}^T + \sigma^2\mathbf{I})^{-1}$ with $\tilde{\mathbf{S}}_{1j}$ and $\mathbf{P}_{1j}$ obtained, respectively, from the matrix $\tilde{\mathbf{S}}_1$ and $\mathbf{P}_1$ by removing the $j$th column. We then have

$$\frac{1}{N}\operatorname{tr}\mathbf{P}_1^2\tilde{\mathbf{S}}_1^T(\tilde{\mathbf{S}}_1\mathbf{P}_1\tilde{\mathbf{S}}_1^T + \sigma^2\mathbf{I})^{-1}\tilde{\mathbf{S}}_1$$

$$(2.5) \qquad = \frac{1}{N}\sum_{j=2}^{K} p_j^2 \breve{\mathbf{s}}_j^T \mathbf{B}_1^{-1} \breve{\mathbf{s}}_j$$

$$= \frac{1}{N}\sum_{j=2}^{K} p_j - \frac{1}{N}\sum_{j=2}^{K} \frac{p_j}{1 + p_j \breve{\mathbf{s}}_j^T B_{1j}^{-1} \breve{\mathbf{s}}_j}.$$

Hence it follows that

$$E\left|\frac{1}{N}\sum_{j=1}^{K}\frac{p_j}{1+p_j\breve{\mathbf{s}}_j^T B_{1j}^{-1}\breve{\mathbf{s}}_j} - \frac{1}{N}\sum_{j=1}^{K}\frac{p_j}{1+p_j\beta}\right|$$

$$\leq M\frac{1}{N}\sum_{j=1}^{K} E\left|\breve{\mathbf{s}}_j^T B_{1j}^{-1}\breve{\mathbf{s}}_j - \frac{1}{N}\operatorname{tr}\mathbf{B}_{1j}^{-1}\right|$$

$$+ M\frac{1}{N}\sum_{j=1}^{K} E\left|\frac{1}{N}\mathbf{B}_{1j}^{-1} - \frac{1}{N}\operatorname{tr}\mathbf{B}_1^{-1}\right|$$

$$(2.6) \qquad + M\frac{K}{N}E\left|\frac{1}{N}\operatorname{tr}\mathbf{B}_1^{-1} - b\right|$$

$$\leq M\frac{cK}{N^{3/2}} + M\frac{1}{N^2}\sum_{j=1}^{K} E\left|\frac{p_j\breve{\mathbf{s}}_j^T \mathbf{B}_{1j}^{-2}\breve{\mathbf{s}}_j}{1+p_j\breve{\mathbf{s}}_j^T \mathbf{B}_{1j}^{-1}\breve{\mathbf{s}}_j}\right| + M\frac{K}{N}E\left|\frac{1}{N}\operatorname{tr}\mathbf{B}_1^{-1} - b\right|$$

$$\to 0 \qquad \text{as } N \to \infty.$$

In the last step, we also use $N^{-1}\operatorname{tr}\mathbf{B}_1^{-1} \xrightarrow{i.p.} b$ and the uniform integrability of $N^{-1}\operatorname{tr}\mathbf{B}_1^{-1}$.

From (2.3)–(2.6) we have

$$(2.7) \qquad \hat{\mathbf{s}}_1^T\hat{\mathbf{S}}_1^T(\hat{\mathbf{S}}_1\hat{\mathbf{S}}_1^T + \sigma^2\mathbf{I})^{-1}\hat{\mathbf{S}}_1\hat{\mathbf{s}}_1 \xrightarrow{i.p.} c\int x\,dH(x) - c\int\frac{x\,dH(x)}{1+xb}.$$

Thus combining (2.3) and (2.7) one can get

$$(2.8) \qquad a_{11} \xrightarrow{i.p.} \left(\sigma^2 + c\int\frac{x\,dH(x)}{1+xb}\right)^{-1} = b.$$

Since $a_{11}$ is bounded by $1/\sigma^2$, $a_{11}$ is then uniformly integrable and so

$$E\left|\frac{1}{N}\sum_{j=1}^{N}a_{jj}^2 - b^2\right| \leq \frac{1}{N}\sum_{j=1}^{N}E|a_{jj}^2 - b^2| = E|a_{11}^2 - b^2| \to 0,$$

as $N \to \infty$. Thus the proof of the above lemma is complete. $\square$

We proceed with the proof of Theorem 1.1. Let $\mathbf{S}_{1k_1}$ and $\mathbf{S}_{1k_1k_2}$ be the matrices obtained from the matrix $\mathbf{S}_1$ by removing the $k_1$th column, the $k_1$th and $k_2$th columns, respectively and let $\mathbf{S}_{1\hat{m}}$ be the matrix obtained from $\mathbf{S}_1$ by deleting the first $m-1$ columns. The matrices $\mathbf{P}_{1k_1}, \mathbf{P}_{1k_1k_2}$ and $\mathbf{P}_{1\hat{m}}$ are defined similarly. Define $\mathbf{A}_{1k_1}^{-1} = (\mathbf{S}_{1k_1}\mathbf{P}_{1k_1}\mathbf{S}_{1k_1}^T + \sigma^2\mathbf{I})^{-1}$, $\mathbf{A}_{1k_1k_2}^{-1} = (\mathbf{S}_{1k_1k_2}\mathbf{P}_{1k_1k_2}\mathbf{S}_{1k_1k_2}^T + \sigma^2\mathbf{I})^{-1}$ and $\mathbf{A}_{1\hat{m}}^{-1} = (\mathbf{S}_{1\hat{m}}\mathbf{P}_{1\hat{m}}\mathbf{S}_{1\hat{m}}^T + \sigma^2\mathbf{I})^{-1}$. Furthermore, all analogues such as $\mathbf{A}_{1k_1k_2k_3}^{-1}$ required in the following derivation have similar meanings.

Write

$$p_1 \mathbf{s}_1^T \mathbf{A}_1^{-1} \mathbf{s}_1$$

$$= p_1 \mathbf{s}_1^T \mathbf{A}_{1\hat{m}}^{-1} \mathbf{s}_1 - \sum_{k_1=2}^{m} \frac{p_1 \mathbf{s}_1^T \mathbf{A}_{1\hat{m}}^{-1} p_{k_1} \mathbf{s}_{k_1} \mathbf{s}_{k_1}^T \mathbf{A}_{1k_1}^{-1} \mathbf{s}_1}{1 + p_{k_1} \mathbf{s}_{k_1}^T \mathbf{A}_{1k_1}^{-1} \mathbf{s}_{k_1}}$$

$$+ \sum_{k_1=2}^{m}\sum_{k_2 \neq k_1}^{m} \frac{p_1 \mathbf{s}_1^T \mathbf{A}_{1\hat{m}}^{-1} p_{k_1} \mathbf{s}_{k_1} \mathbf{s}_{k_1}^T \mathbf{A}_{1\hat{m}}^{-1} p_{k_2} \mathbf{s}_{k_2} \mathbf{s}_{k_2}^T \mathbf{A}_{1k_1k_2}^{-1} \mathbf{s}_1}{(1 + p_{k_1} \mathbf{s}_{k_1}^T \mathbf{A}_{1k_1}^{-1} \mathbf{s}_{k_1})(1 + p_{k_2} \mathbf{s}_{k_2}^T \mathbf{A}_{1k_1k_2}^{-1} \mathbf{s}_{k_2})}$$

(2.9) $$= p_1 \mathbf{s}_1^T \mathbf{A}_{1\hat{m}}^{-1} \mathbf{s}_1 - \sum_{k_1=2}^{m} \frac{p_1 \mathbf{s}_1^T \mathbf{A}_{1\hat{m}}^{-1} p_{k_1} \mathbf{s}_{k_1} \mathbf{s}_{k_1}^T \mathbf{A}_{1\hat{m}}^{-1} \mathbf{s}_1}{1 + p_{k_1} \mathbf{s}_{k_1}^T \mathbf{A}_{1k_1}^{-1} \mathbf{s}_{k_1}}$$

$$+ \sum_{k_1=2}^{m}\sum_{k_2 \neq k_1}^{m} \frac{p_1 \mathbf{s}_1^T \mathbf{A}_{1\hat{m}}^{-1} p_{k_1} \mathbf{s}_{k_1} \mathbf{s}_{k_1}^T \mathbf{A}_{1\hat{m}}^{-1} p_{k_2} \mathbf{s}_{k_2} \mathbf{s}_{k_2}^T \mathbf{A}_{1\hat{m}}^{-1} \mathbf{s}_1}{(1 + p_{k_1} \mathbf{s}_{k_1}^T \mathbf{A}_{1k_1}^{-1} \mathbf{s}_{k_1})(1 + p_{k_2} \mathbf{s}_{k_2}^T \mathbf{A}_{1k_1k_2}^{-1} \mathbf{s}_{k_2})} + \cdots$$

$$+ (-1)^{m+1}$$

$$\times \sum_{\substack{k_1=2, \\ k_2 \neq k_1, \ldots, \\ k_{m-1} \neq k_{m-2}}}^{m} \frac{p_1 \mathbf{s}_1^T \mathbf{A}_{1\hat{m}}^{-1} p_{k_1} \mathbf{s}_{k_1} \mathbf{s}_{k_1}^T \mathbf{A}_{1\hat{m}}^{-1} p_{k_2} \mathbf{s}_{k_2} \cdots p_{k_{m-1}} \mathbf{s}_{k_{m-1}}^T \mathbf{A}_{1\hat{m}}^{-1} \mathbf{s}_1}{(1 + p_{k_1} \mathbf{s}_{k_1}^T \mathbf{A}_{1k_1}^{-1} \mathbf{s}_{k_1}) \cdots (1 + p_{k_{m-1}} \mathbf{s}_{k_{m-1}}^T \mathbf{A}_{1\hat{m}}^{-1} \mathbf{s}_{k_{m-1}})},$$

where the subscripts $k_1, \ldots, k_{m-1}$ are larger than 1.

For any $i \neq j$ $(i, j = 1, \ldots, m)$ we have

$$E(\mathbf{s}_i^T \mathbf{A}_{1\hat{m}}^{-1} \mathbf{s}_j)^2 = E \operatorname{tr} \mathbf{A}_{1\hat{m}}^{-1} \mathbf{s}_j \mathbf{s}_j^T \mathbf{A}_{1\hat{m}}^{-1} \mathbf{s}_i \mathbf{s}_i^T = E \frac{1}{N^2} \operatorname{tr} \mathbf{A}_{1\hat{m}}^{-2}$$

which implies

$$N^{1/4} \mathbf{s}_i^T A_{1\hat{m}}^{-1} \mathbf{s}_j \xrightarrow{i.p.} 0$$

and then

(2.10) $$\sqrt{N}\mathbf{s}_1^T \mathbf{A}_{1\hat{m}}^{-1}\mathbf{s}_{k_1}\mathbf{s}_{k_1}^T \mathbf{A}_{1\hat{m}}^{-1}\mathbf{s}_{k_2} \xrightarrow{i.p.} 0.$$

Hence, from (2.9) and (2.10) we have

(2.11) $$\sqrt{N}(p_1\mathbf{s}_1^T \mathbf{A}_1^{-1}\mathbf{s}_1 - \sqrt{N}p_1\mathbf{s}_1^T \mathbf{A}_{1\hat{m}}^{-1}\mathbf{s}_1) \xrightarrow{i.p.} 0.$$

Similarly, for any $k = 1, \ldots, m$ one can find

$$\sqrt{N}(p_k\mathbf{s}_k^T \mathbf{A}_k^{-1}\mathbf{s}_k - \sqrt{N}p_k\mathbf{s}_k^T \mathbf{A}_{1\hat{m}}^{-1}\mathbf{s}_k) \xrightarrow{i.p.} 0.$$

It thus suffices to consider the asymptotic distribution for the linear combination of $\sqrt{N}(p_k\mathbf{s}_k^T \mathbf{A}_{1\hat{m}}^{-1}\mathbf{s}_k - p_k b_N)$, $k = 1, \ldots, m$. To this end, it can be seen that $\{p_k\mathbf{s}_k^T \mathbf{A}_{1\hat{m}}^{-1}\mathbf{s}_k, k = 1, \ldots, m\}$ are independent when the matrix $\mathbf{A}_{1\hat{m}}^{-1}$ is given and hence it suffices to consider only one of $\{p_k\mathbf{s}_k^T \mathbf{A}_{1\hat{m}}^{-1}\mathbf{s}_k, k = 1, \ldots, m\}$ when the matrix $\mathbf{A}_{1\hat{m}}^{-1}$ is given.

By Lemma 2.1 and the Jensen inequality, it is easy to verify that

$$\lim_{N\to\infty} \sum_j^N a_{jj}^2 \Big/ \sum_{i,j} a_{ij}^2 < 1,$$

where $\mathbf{A}_{1\hat{m}}^{-1} = (a_{ij})$. Hence for any $k = 1, \ldots, m$,

$$\sqrt{N}\left(p_k\mathbf{s}_k^T \mathbf{A}_{1\hat{m}}^{-1}\mathbf{s}_k - p_k\frac{\operatorname{tr} \mathbf{A}_{1\hat{m}}^{-1}}{N}\right)$$
$$\xrightarrow{D} N\left(0, \left(2\int \frac{dF^{c,H}(x)}{(x+\sigma^2)^2} + \left(\int \frac{dF^{c,H}(x)}{x+\sigma^2}\right)^2 (Ev_{11}^4 - 3)\right)p_k^2\right)$$

by Theorem 1.1 in [5] when $\mathbf{A}_{1\hat{m}}^{-1}$ is given. Here the asymptotic variance can be computed using formula (4.23). From result (1) of Theorem 1.1 of [3] it can be concluded that

$$\sqrt{N}\left(p_k\frac{\operatorname{tr} \mathbf{A}_{1\hat{m}}^{-1}}{N} - p_k b_N\right) \xrightarrow{i.p.} 0.$$

Thus we are done by the Fubini theorem and the Cramér–Wold device.

**3. Proof of Theorem 1.2.** Let $z = u + iv$, $v > 0$. Recall that the Stieltjes transform is defined for any distribution function $F$ as

$$m_F(z) = \int \frac{1}{x-z} dF(x), \qquad z \in \mathbb{C}^+ \equiv \{z \in \mathbb{C}, \Im z > 0\}.$$

Hence, the Stieltjes transform of $G_N(x)$ is

(3.1) $$m_{G_N(x)}(z) = \frac{1}{K}\sum_{k=1}^K \frac{1}{\beta_k - z}$$

and it suffices to consider $m_{G_N(x)}(z)$.

First we obtain a decomposition as follows:

$$(3.2) \quad m_{G_N(x)}(z) - \int \frac{dH(x)}{bx - z} = V_1 + V_2 + V_3 + V_4,$$

where

$$V_1 = m_{G_N(x)}(z) - \frac{1}{K} \sum_{k=1}^{K} \frac{1}{(p_k/N) \operatorname{tr} \mathbf{A}_k^{-1} - z},$$

$$V_2 = \frac{1}{K} \sum_{k=1}^{K} \left( \frac{1}{(p_k/N) \operatorname{tr} \mathbf{A}_k^{-1} - z} - \frac{1}{(p_k/N) \operatorname{tr} \mathbf{A}^{-1} - z} \right),$$

$$V_3 = \frac{1}{K} \sum_{k=1}^{K} \left( \frac{1}{(p_k/N) \operatorname{tr} \mathbf{A}^{-1} - z} - \frac{1}{p_k b - z} \right),$$

$$V_4 = \frac{1}{K} \sum_{k=1}^{K} \frac{1}{p_k b - z} - \int \frac{dH(x)}{xb - z}.$$

It is straightforward to verify that

$$(3.3) \quad \left| \frac{1}{\beta_k - z} \right| \leq \frac{1}{v}, \quad \left| \frac{1}{p_k b - z} \right| \leq \frac{1}{v};$$

then we have

$$(3.4) \quad E|V_1| \leq \frac{1}{v^2 K} \sum_{k=1}^{K} p_k E \left| \mathbf{s}_k^T \mathbf{A}_k^{-1} \mathbf{s}_k - \frac{1}{N} \operatorname{tr} \mathbf{A}_k^{-1} \right| \leq \frac{M}{N} \to 0,$$

as $N \to \infty$. Similarly, by (3.3) one can find

$$(3.5) \quad E|V_2| \leq \frac{1}{v^2 NK} \sum_{k=1}^{K} E \frac{p_k^2 \mathbf{s}_k^T \mathbf{A}_k^{-2} \mathbf{s}_k}{1 + p_k \mathbf{s}_k^T \mathbf{A}_k^{-1} \mathbf{s}_k} \leq \frac{M}{N} \to 0,$$

as $N \to \infty$. From the uniform integrability of the random variable $N^{-1} \operatorname{tr} \mathbf{A}^{-1}$, one can obtain

$$(3.6) \quad E|V_3| \leq \frac{1}{v^2 K} \sum_{k=1}^{K} p_k E \left| \frac{1}{N} \operatorname{tr} \mathbf{A}^{-1} - b \right| \to 0.$$

It is obvious that $|V_4|$ converges to zero. This, together with (3.1)–(3.6), implies Theorem 1.2 and thus we are done.

**4. Proof of Theorem 1.3.** We begin the proof of this theorem with the replacement of the entries of $\mathbf{S}$ by truncated and centralized variables. Since $Ev_{11}^4 < \infty$, we have $\varepsilon^{-4} Ev_{11}^4 I(|v_{11}| > \varepsilon \sqrt{N}) \to 0$ for any $\varepsilon > 0$. Thus a positive sequence $\varepsilon_N$ converging to zero can be selected so that

(4.1) $$\varepsilon_N^{-4} Ev_{11}^4 I(|v_{11}| > \varepsilon_N \sqrt{N}) \to 0.$$

Define $\hat{v}_{ij} = v_{ij} I(|v_{ij}| \leq \varepsilon_N \sqrt{N})$ and $\bar{v}_{ij} = \hat{v}_{ij} - E\hat{v}_{ij}$, $i = 1, \ldots, N$, $j = 1, \ldots, K$. The corresponding matrices and vectors are denoted by $\hat{\mathbf{s}}_k, \bar{\mathbf{s}}_k, \hat{\mathbf{S}}_k$ and $\bar{\mathbf{S}}_k$, $k = 1, \ldots, K$, the elements of which are $\hat{v}_{ij}$ or $\bar{v}_{ij}$ instead of $v_{ij}$.

Let
$$\hat{\mathbf{A}}_k = \hat{\mathbf{S}}_k \mathbf{P}_k \hat{\mathbf{S}}_k^T + \sigma^2 \mathbf{I}, \qquad \bar{\mathbf{A}}_k = \bar{\mathbf{S}}_k \mathbf{P}_k \bar{\mathbf{S}}_k^T + \sigma^2 \mathbf{I}.$$

It follows from (4.1) that

(4.2) $$P\left(\sum_{k=1}^K p_k \mathbf{s}_k^T \mathbf{A}_k^{-1} \mathbf{s}_k \neq \sum_{k=1}^K p_k \hat{\mathbf{s}}_k^T \hat{\mathbf{A}}_k^{-1} \hat{\mathbf{s}}_k\right)$$
$$\leq P\left(\bigcup_{i,j}(|v_{ij}| \geq \varepsilon_N \sqrt{N})\right)$$
$$\leq NKP(|v_{11}| \geq \varepsilon_N \sqrt{N}) \to 0,$$

as $N \to \infty$.

Observe that

(4.3) $$|\hat{\mathbf{s}}_k^T \hat{\mathbf{A}}_k^{-1} \hat{\mathbf{s}}_k - \bar{\mathbf{s}}_k^T \hat{\mathbf{A}}_k^{-1} \bar{\mathbf{s}}_k)| \leq (E\hat{\mathbf{s}}_k^T) \hat{\mathbf{A}}_k^{-1} E\hat{\mathbf{s}}_k + 2|\hat{\mathbf{s}}_k^T \hat{\mathbf{A}}_k^{-1} E\hat{\mathbf{s}}_k|$$
$$\leq \frac{1}{\sigma^2} \|E\hat{\mathbf{s}}_k\|^2 + 2|\hat{\mathbf{s}}_k^T \hat{\mathbf{A}}_k^{-1} E\hat{\mathbf{s}}_k|.$$

Concerning the first item on the right, we have

(4.4) $$\|E\hat{\mathbf{s}}_k\|^2 = \|E\hat{\mathbf{s}}_1\|^2 = o(N^{-3}),$$

where we use the fact that

(4.5) $$E\hat{v}_{11} = o(N^{-3/2}).$$

For the second item, by (4.4) one can find

(4.6) $$E(\bar{\mathbf{s}}_k^T \hat{\mathbf{A}}_k^{-1} E\hat{\mathbf{s}}_k)^2 \leq \frac{1}{\sigma^4} \|E\hat{\mathbf{s}}_k\|^2 E\|\bar{\mathbf{s}}_k\|^2$$
$$= \frac{1}{\sigma^4} \|E\hat{\mathbf{s}}_1\|^2 E\|\bar{\mathbf{s}}_1\|^2 = o(N^{-3}).$$

Combining (4.3)–(4.6), one can conclude that

$$E\left|\sum_k^K (\hat{\mathbf{s}}_k^T \hat{\mathbf{A}}_k^{-1} \hat{\mathbf{s}}_k - \bar{\mathbf{s}}_k^T \hat{\mathbf{A}}_k^{-1} \bar{\mathbf{s}}_k)\right|$$

(4.7)
$$\leq \frac{3K}{\sigma^2} \|E\hat{\mathbf{s}}_1\|^2 + 2\sum_k^K (E(\hat{\mathbf{s}}_k^T \hat{\mathbf{A}}_k^{-1} E\hat{\mathbf{s}}_k)^2)^{1/2} \longrightarrow 0.$$

Next we will show that

(4.8)
$$\sum_k^K (\bar{\mathbf{s}}_k^T \hat{\mathbf{A}}_k^{-1} \bar{\mathbf{s}}_k - \bar{\mathbf{s}}_k^T \bar{\mathbf{A}}_k^{-1} \bar{\mathbf{s}}_k) \xrightarrow{i.p.} 0.$$

By matrix inverse formula $\mathbf{A}^{-1} - \mathbf{B}^{-1} = \mathbf{B}^{-1}(\mathbf{B} - \mathbf{A})\mathbf{A}^{-1}$, we have

$$|\bar{\mathbf{s}}_k^T \hat{\mathbf{A}}_k^{-1} \bar{\mathbf{s}}_k - \bar{\mathbf{s}}_k^T \bar{\mathbf{A}}_k^{-1} \bar{\mathbf{s}}_k|$$
$$= |\bar{\mathbf{s}}_k^T \bar{\mathbf{A}}_k^{-1} [(E\hat{\mathbf{S}}_k)\mathbf{P}_k E(\hat{\mathbf{S}}_k^T)$$
$$\quad - (E\hat{\mathbf{S}}_k)\mathbf{P}_k \hat{\mathbf{S}}_k^T - \hat{\mathbf{S}}_k \mathbf{P}_k (E\hat{\mathbf{S}}_k^T)]\hat{\mathbf{A}}_k^{-1} \bar{\mathbf{s}}_k|$$

(4.9)
$$\leq \frac{1}{\sigma^4} \|\bar{\mathbf{s}}_k\|^2 \|E\hat{\mathbf{S}}_k\|^2 \|\mathbf{P}_k\|$$
$$\quad + |\bar{\mathbf{s}}_k^T \bar{\mathbf{A}}_k^{-1} [(E\hat{\mathbf{S}}_k)\mathbf{P}_k \hat{\mathbf{S}}_k^T + \hat{\mathbf{S}}_k \mathbf{P}_k (E\hat{\mathbf{S}}_k^T)]\hat{\mathbf{A}}_k^{-1} \bar{\mathbf{s}}_k|$$
$$\leq \frac{MN}{\sigma^4} \|\bar{\mathbf{s}}_k\|^2 (E\hat{v}_{11})^2$$
$$\quad + |\bar{\mathbf{s}}_k^T \bar{\mathbf{A}}_k^{-1} [(E\hat{\mathbf{S}}_k)\mathbf{P}_k \hat{\mathbf{S}}_k^T + \hat{\mathbf{S}}_k \mathbf{P}_k E(\hat{\mathbf{S}}_k^T)]\hat{\mathbf{A}}_k^{-1} \bar{\mathbf{s}}_k|.$$

By Lemma 2.7 of [2] one can obtain

$$E\left(\bar{\mathbf{s}}_k^T \bar{\mathbf{A}}_k^{-1} \hat{\mathbf{S}}_k \mathbf{P}_k (E\hat{\mathbf{S}}_k^T) \hat{\mathbf{A}}_k^{-1} \bar{\mathbf{s}}_k - \frac{1}{N} \operatorname{tr} \bar{\mathbf{A}}_k^{-1} \hat{\mathbf{S}}_k \mathbf{P}_k (E\hat{\mathbf{S}}_k^T) \hat{\mathbf{A}}_k^{-1}\right)^2$$
$$\leq \frac{NM^2}{N^2 \sigma^8} (E\hat{v}_{11})^2 E \operatorname{tr} \hat{\mathbf{S}}_k \hat{\mathbf{S}}_k^T$$
$$\leq \frac{NM^2}{N^2 \sigma^8} (E\hat{v}_{11})^2 E \operatorname{tr} \hat{\mathbf{S}}_1 \hat{\mathbf{S}}_1^T,$$

which implies

(4.10) $$\sum_{k=1}^K \bar{\mathbf{s}}_k^T \bar{\mathbf{A}}_k^{-1} \hat{\mathbf{S}}_k \mathbf{P}_k (E\hat{\mathbf{S}}_k^T) \hat{\mathbf{A}}_k^{-1} \bar{\mathbf{s}}_k - \frac{1}{N} \operatorname{tr} \bar{\mathbf{A}}_k^{-1} \hat{\mathbf{S}}_k \mathbf{P}_k (E\hat{\mathbf{S}}_k^T) \hat{\mathbf{A}}_k^{-1} \xrightarrow{i.p.} 0.$$

Here we use (4.5) and

(4.11) $$\frac{1}{N} E \operatorname{tr} \hat{\mathbf{S}}_1 \hat{\mathbf{S}}_1^T = \frac{1}{N} E \sum_{j=2}^K \hat{\mathbf{s}}_j^T \hat{\mathbf{s}}_j \leq \frac{K^2}{N^2} E v_{11}^2.$$

On the other hand we have

$$\left|\frac{1}{N}\operatorname{tr}\bar{\mathbf{A}}_k^{-1}\hat{\mathbf{S}}_k\mathbf{P}_k E(\hat{\mathbf{S}}_k^T)\hat{\mathbf{A}}_k^{-1}\right| = \frac{|E\hat{v}_{11}|}{N^{3/2}}\left|\sum_{j=1,j\neq k} p_j\mathbf{e}^T\hat{\mathbf{A}}_k^{-1}\bar{\mathbf{A}}_k^{-1}\hat{\mathbf{s}}_j\right|$$

$$\leq \frac{M|E\hat{v}_{11}|}{\sigma^4 N}\sum_{j=1}^K \|\hat{\mathbf{s}}_j\|,$$

where $\mathbf{e} = (1,\ldots,1)^T$. This, together with (4.5), gives

$$E\left|\frac{1}{N}\sum_{k=1}^K \operatorname{tr}\bar{\mathbf{A}}_k^{-1}\hat{\mathbf{S}}_k\mathbf{P}_k(E\hat{\mathbf{S}}_k^T)\hat{\mathbf{A}}_k^{-1}\right| \leq \frac{KM|E\hat{v}_{11}|}{\sigma^4 N}\sum_{j=1}^K (E\|\hat{\mathbf{s}}_j\|^2)^{1/2} \xrightarrow{i.p.} 0,$$

and by combining (4.10) one can then find

(4.12) $$\sum_{k=1}^K \bar{\mathbf{s}}_k^T \bar{\mathbf{A}}_k^{-1}\hat{\mathbf{S}}_k\mathbf{P}_k(E\hat{\mathbf{S}}_k^T)\hat{\mathbf{A}}_k^{-1}\bar{\mathbf{s}}_k \xrightarrow{i.p.} 0.$$

Similarly, one can also show that

(4.13) $$\sum_{k=1}^K \bar{\mathbf{s}}_k^T \bar{\mathbf{A}}_k^{-1}(E\hat{\mathbf{S}}_k)\mathbf{P}_k\hat{\mathbf{S}}_k^T\hat{\mathbf{A}}_k^{-1}\bar{\mathbf{s}}_k \xrightarrow{i.p.} 0.$$

Thus (4.8) immediately follows from (4.5), (4.9), (4.12) and (4.13). It is easy to check that

(4.14) $$1 - \operatorname{var}(\bar{v}_{11}) = o(N^{-1}).$$

Applying this and the argument similar to the centralization step, one can then renormalize the underlying random variables. Consequently, it can be assumed that the underlying random variables satisfy

$$Ev_{11} = 0, \qquad Ev_{11}^2 = 1, \qquad |v_{11}| \leq \varepsilon_N\sqrt{N}.$$

In the sequel we still use $v_{ij}$, $\mathbf{s}_k$, $\mathbf{S}_k$ and $\mathbf{A}_k$ instead of $\bar{v}_{11}$, $\bar{\mathbf{s}}_k$, $\bar{\mathbf{S}}_k$ and $\bar{\mathbf{A}}_k$ to simplify the notation.

Write $\hat{s}_k = \mathbf{s}_k^T\mathbf{A}_k^{-1}\mathbf{s}_k$ and

$$\sum_{k=1}^K p_k\hat{s}_k = \sum_{k=1}^K p_k\mathbf{s}_k^T(\mathbf{A}_k^{-1} - \mathbf{A}^{-1})\mathbf{s}_k + \operatorname{tr}\left(\mathbf{A}^{-1}\left(\sum_{k=1}^K p_k\mathbf{s}_k\mathbf{s}_k^T\right)\right)$$

$$= \sum_{k=1}^K \frac{(p_k\hat{s}_k)^2}{1+p_k\hat{s}_k} + N - \sigma^2\operatorname{tr}\mathbf{A}^{-1}.$$

Further, after some simple computations one can find

$$-\sum_{k=1}^K \frac{1}{1+p_k\hat{s}_k} = N - \sigma^2\operatorname{tr}\mathbf{A}^{-1} - K.$$

Applying the formula

$$\text{(4.15)} \quad \frac{1}{1+p_k\hat{s}_k} = \frac{1}{1+b_Np_k} - \frac{p_k\hat{s}_k - b_Np_k}{(1+p_k\hat{s}_k)(1+b_Np_k)}$$

to the above identity we can arrive at

$$\text{(4.16)} \quad \int \frac{dH_N(x)}{(1+xb_N)^2} \sum_{k=1}^{K} p_k(\hat{s}_k - b_N) = U_1 - U_2 - U_3 + U_4 - U_5,$$

where

$$U_1 = \sum_{k=1}^{K} \frac{p_k^2(\hat{s}_k - b_N)^2}{(1+b_Np_k)^3},$$

$$U_2 = \sigma^2(\operatorname{tr} A^{-1} - Nb_N),$$

$$U_3 = \sum_{k=1}^{K} \frac{p_k^3(\hat{s}_k - b_N)^3}{(1+p_kb_N)^3(1+p_k\hat{s}_k)},$$

$$U_4 = N(1 - \sigma^2 b_N) - K + K \int \frac{dH_N(x)}{1+xb_N}$$

and

$$U_5 = \sum_{k=1}^{K} p_k(\hat{s}_k - b_N)\left(\frac{1}{(1+p_kb_N)^2} - \int \frac{dH_N(x)}{(1+xb_N)^2}\right).$$

As will be seen, the contributions from $U_3$, $U_4$ and $U_5$ can be ignored and the main terms are $U_1$ and $U_2$.

It is easy to see that $b_N$ satisfies

$$\text{(4.17)} \quad \frac{1}{b_N} = \sigma^2 + c_N \int \frac{x\,dH_N(x)}{1+xb_N},$$

that is, $U_4 = 0$.

For the term $U_3$ we have

$$|U_3| \leq M \sum_{k=1}^{K} |\hat{s}_k - b_N|^3 \leq M(U_{31} + U_{32} + U_{33}),$$

where

$$U_{31} = \sum_{k=1}^{K} \left|\hat{s}_k - \frac{1}{N}\operatorname{tr} \mathbf{A}_k^{-1}\right|^3,$$

$$U_{32} = \sum_{k=1}^{K} \left|\frac{1}{N}\operatorname{tr} \mathbf{A}_k^{-1} - \frac{1}{N}\operatorname{tr} \mathbf{A}^{-1}\right|^3$$

and
$$U_{33} = K \left| \frac{1}{N} \operatorname{tr} \mathbf{A}^{-1} - b_N \right|^3.$$

From Lemma 2.7 of [2] one can find
$$EU_{31} \leq \sum_{k=1}^{K} \left( \frac{M}{N^3} Ev_{11}^4 E(\operatorname{tr} \mathbf{A}_k^{-2})^{3/2} + \frac{M}{N^3} Ev_{11}^6 E \operatorname{tr} \mathbf{A}_k^{-3} \right) = O(\varepsilon_N).$$

For the term $U_{32}$
$$U_{32} \leq \frac{1}{N^3} \sum_{k=1}^{K} \left| \frac{p_k \mathbf{s}_k^T \mathbf{A}_k^{-2} \mathbf{s}_k}{1 + p_k \hat{s}_k} \right|^3 = O(N^{-2}),$$

where we use the fact

(4.18) $$\left| \frac{p_k \mathbf{s}_k^T \mathbf{A}_k^{-2} \mathbf{s}_k}{1 + p_k \hat{s}_k} \right| \leq \frac{1}{\sigma^2}.$$

In the sequel, we will not mention it again whenever (4.18) is used. By Theorem 1 of [3] we have
$$U_{33} \xrightarrow{i.p.} 0.$$

From the above argument it can be concluded that

(4.19) $$U_3 \xrightarrow{i.p.} 0.$$

We now analyze the term $U_1$ by computing its variance:

(4.20) $$E \left( \sum_{k=1}^{K} \frac{p_k^2 (\hat{s}_k - b_N)^2}{(1 + b_N p_k)^3} - \sum_{k=1}^{K} \frac{p_k^2 E(\hat{s}_k - b_N)^2}{(1 + b_N p_k)^3} \right)^2 = U_{11} + U_{12},$$

where
$$U_{11} = \sum_{k=1}^{K} E \left( \frac{p_k^2 (\hat{s}_k - b_N)^2}{(1 + b_N p_k)^3} - \frac{p_k^2 E(\hat{s}_k - b_N)^2}{(1 + b_N p_k)^3} \right)^2$$

and
$$U_{12} = \sum_{k_1 \neq k_2}^{K} E \left[ \left( \frac{p_{k_1}^2 (\hat{s}_{k_1} - b_N)^2}{(1 + b_N p_{k_1})^3} - \frac{p_{k_1}^2 E(\hat{s}_{k_1} - b_N)^2}{(1 + b_N p_{k_1})^3} \right) \right.$$
$$\left. \times \left( \frac{p_{k_2}^2 (\hat{s}_{k_2} - b_N)^2}{(1 + b_N p_{k_2})^3} - \frac{p_{k_2}^2 E(\hat{s}_{k_2} - b_N)^2}{(1 + b_N p_{k_2})^3} \right) \right].$$

Similarly to the argument of (4.19), one can get

$$U_{11} \leq M \sum_{k=1}^{K} E(\hat{s}_k - b_N)^4$$

$$\leq M \sum_{k=1}^{K} \left[ E\left(\hat{s}_k - \frac{1}{N} \operatorname{tr} \mathbf{A}_k^{-1}\right)^4 + KE\left(\frac{1}{N} \operatorname{tr} \mathbf{A}_k^{-1} - \frac{1}{N} \operatorname{tr} \mathbf{A}^{-1}\right)^4 \right]$$

(4.21)
$$+ MKE\left(\frac{1}{N} \operatorname{tr} \mathbf{A}^{-1} - b_N\right)^4$$

$$= O(\varepsilon_N),$$

as $N \to \infty$. Indeed, in the last step we also use the fact that

$$KE\left(\frac{1}{N} \operatorname{tr} \mathbf{A}^{-1} - b_N\right)^4$$

$$\leq \left(1/\sigma^2 + b_N\right)^2 \frac{K}{N^2} E(\operatorname{tr} \mathbf{A}^{-1} - Nb_N)^2 \to 0.$$

To evaluate the term $U_{12}$, we need to decompose it further as shown below. The strategy is to split $\mathbf{A}_{k_1}^{-1}$ into the sum of $\mathbf{A}_{k_1 k_2}^{-1}$ and

$$-\frac{\mathbf{A}_{k_1 k_2}^{-1} p_{k_2} \mathbf{s}_{k_2} \mathbf{s}_{k_2}^T \mathbf{A}_{k_1 k_2}^{-1}}{1 + p_{k_2} \mathbf{s}_{k_2}^T \mathbf{A}_{k_1 k_2}^{-1} \mathbf{s}_{k_2}},$$

so does for $\mathbf{A}_{k_2}^{-1}$. Thus one can find

(4.22) $$U_{12} = U_{121} + U_{122} + U_{123} + U_{124} + \cdots + U_{129},$$

where

$$U_{121} = \sum_{k_1 \neq k_2}^{K} E(p_{k_1 k_2} \beta_{k_1 k_2} \beta_{k_2 k_1}), \qquad U_{122} = \sum_{k_1 \neq k_2}^{K} E(p_{k_1 k_2} \zeta_{k_1 k_2} \beta_{k_2 k_1}),$$

$$U_{123} = \sum_{k_1 \neq k_2}^{K} E(p_{k_1 k_2} \beta_{k_1 k_2} \zeta_{k_2 k_1}), \qquad U_{124} = \sum_{k_1 \neq k_2}^{K} E(p_{k_1 k_2} \zeta_{k_1 k_2} \zeta_{k_2 k_1}),$$

$$U_{125} = -2 \sum_{k_1 \neq k_2}^{K} E(p_{k_1 k_2} \zeta_{k_1 k_2} \alpha_{k_2 k_1}), \qquad U_{126} = -2 \sum_{k_1 \neq k_2}^{K} E(p_{k_1 k_2} \beta_{k_1 k_2} \alpha_{k_2 k_1}),$$

$$U_{127} = -2 \sum_{k_1 \neq k_2}^{K} E(p_{k_1 k_2} \alpha_{k_1 k_2} \zeta_{k_2 k_1}), \qquad U_{128} = -2 \sum_{k_1 \neq k_2}^{K} E(p_{k_1 k_2} \alpha_{k_1 k_2} \beta_{k_2 k_1}),$$

$$U_{129} = 4 \sum_{k_1 \neq k_2}^{K} E(p_{k_1 k_2} \alpha_{k_1 k_2} \alpha_{k_2 k_1}),$$

with

$$p_{ij} = \frac{1}{(1+b_N p_i)^3 (1+b_N p_j)^3},$$

$$\beta_{ij} = p_i^2 (\mathbf{s}_i^T \mathbf{A}_{ij}^{-1} \mathbf{s}_i - b_N)^2 - p_i^2 E(\mathbf{s}_i^T \mathbf{A}_{ij}^{-1} \mathbf{s}_i - b_N)^2,$$

$$\zeta_{ij} = \frac{\gamma_{ij}^2}{(1+p_j \mathbf{s}_j^T \mathbf{A}_{ij}^{-1} \mathbf{s}_j)^2} - E \frac{\gamma_{ij}^2}{(1+p_j \mathbf{s}_j^T \mathbf{A}_{ij}^{-1} \mathbf{s}_j)^2},$$

$$\gamma_{ij} = p_i p_j \mathbf{s}_i^T \mathbf{A}_{ij}^{-1} \mathbf{s}_j \mathbf{s}_j^T \mathbf{A}_{ij}^{-1} \mathbf{s}_i,$$

$$\alpha_{ij} = \frac{p_i \gamma_{ij} (\mathbf{s}_i^T \mathbf{A}_{ij}^{-1} \mathbf{s}_i - b_N)}{1 + p_j \mathbf{s}_j^T \mathbf{A}_{ij}^{-1} \mathbf{s}_j} - E \frac{p_i \gamma_{ij} (\mathbf{s}_i^T \mathbf{A}_{ij}^{-1} \mathbf{s}_i - b_N)}{1 + p_j \mathbf{s}_j^T \mathbf{A}_{ij}^{-1} \mathbf{s}_j}.$$

Also, set

$$\alpha_{k_1} = \frac{1}{1+p_{k_1} \mathbf{s}_{k_1}^T \mathbf{A}_{k_1 k_2}^{-1} \mathbf{s}_{k_1}}, \qquad \alpha_{k_2} = \frac{1}{1+p_{k_2} \mathbf{s}_{k_2}^T \mathbf{A}_{k_1 k_2}^{-1} \mathbf{s}_{k_2}}.$$

As will be seen, each of $\beta_{ij}, \zeta_{ij}$ and $\alpha_{ij}$ converges to zero in some way and the convergence rate is needed to attain our aim. In the subsequent paragraphs we show that each term $U_{12j}, j=1,\ldots,9$, converges to zero.

Consider the term $U_{121}$ first. It is straightforward to verify that

(4.23)
$$E \left| \mathbf{s}_k^T \mathbf{B}_k \mathbf{s}_k - \frac{1}{N} \operatorname{tr} \mathbf{B}_k \right|^2 = \frac{1}{N^2}(Ev_{11}^4 - 3) \sum_{j=1}^{N} E(b_{jj}^{(k)})^2 + \frac{2}{N^2} E \operatorname{tr} \mathbf{B}_k \mathbf{B}_k^T,$$

where $\mathbf{B}_k = (b_{j_1 j_2}^{(k)})$ is any symmetric matrix independent of $\mathbf{s}_k$. It follows that

$$E[((\mathbf{s}_{k_1}^T \mathbf{A}_{k_1 k_2}^{-1} \mathbf{s}_{k_1} - b_N)^2 - E(\mathbf{s}_{k_1}^T \mathbf{A}_{k_1 k_2}^{-1} \mathbf{s}_{k_1} - b_N)^2) | \mathbf{A}_{k_1 k_2}^{-1}]$$
$$= \frac{2}{N^2} \operatorname{tr}(\mathbf{A}_{k_1 k_2}^{-2} - E \mathbf{A}_{k_1 k_2}^{-2}) + \left( \frac{1}{N} \operatorname{tr} \mathbf{A}_{k_1 k_2}^{-1} - b_N \right)^2 - E \left( \frac{1}{N} \operatorname{tr} \mathbf{A}_{k_1 k_2}^{-1} - b_N \right)^2,$$

and then, that

$$U_{121} = E \left[ \sum_{k_1 \neq k_2}^{K} p_{k_1 k_2} E(\beta_{k_1 k_2} | \mathbf{A}_{k_1 k_2}^{-1}) E(\beta_{k_2 k_1} | \mathbf{A}_{k_1 k_2}^{-1}) \right],$$

(4.24)
$$\leq \frac{M}{N^4} \sum_{k_1 \neq k_2}^{K} (E(\operatorname{tr} \mathbf{A}_{k_1 k_2}^{-2} - \operatorname{tr} E \mathbf{A}_{k_1 k_2}^{-2})^2 + E(\operatorname{tr} \mathbf{A}_{k_1 k_2}^{-1} - N b_N)^4 + (E(\operatorname{tr} \mathbf{A}_{k_1 k_2}^{-1} - N b_N)^2)^2).$$

Since the distribution of $\operatorname{tr} \mathbf{A}_{k_1 k_2}^{-2}$ is dependent on different $k_1, k_2$, the difference between $\operatorname{tr} \mathbf{A}_{k_1 k_2}^{-2}$ and $\operatorname{tr} \mathbf{A}^{-2}$ caused by a different $k$ must be eliminated. To this end, by splitting $\mathbf{A}_{k_1 k_2}^{-1}$ into the sum of $\mathbf{A}_{k_1}^{-1}$ and $\xi_{k_1 k_2}$, one can get

$$\begin{aligned}
&E(\operatorname{tr} \mathbf{A}_{k_1 k_2}^{-2} - \operatorname{tr} E\mathbf{A}_{k_1 k_2}^{-2})^2 \\
&\leq ME(\operatorname{tr} \mathbf{A}_{k_1 k_2}^{-1} \xi_{k_1 k_2} - \operatorname{tr} E\mathbf{A}_{k_1 k_2}^{-1} \xi_{k_1 k_2})^2 \\
&\quad + ME(\operatorname{tr} \xi_{k_1 k_2} \mathbf{A}_{k_1}^{-1} - \operatorname{tr} E\xi_{k_1 k_2} \mathbf{A}_{k_1}^{-1})^2 \\
&\quad + ME(\operatorname{tr} \mathbf{A}_{k_1}^{-2} - \operatorname{tr} E\mathbf{A}_{k_1}^{-2})^2 \\
&\leq M + ME(\operatorname{tr} \mathbf{A}_{k_1}^{-2} - \operatorname{tr} E\mathbf{A}_{k_1}^{-2})^2,
\end{aligned} \quad (4.25)$$

where

$$\xi_{k_1 k_2} = \mathbf{A}_{k_1 k_2}^{-1} p_{k_2} \mathbf{s}_{k_2} \mathbf{s}_{k_2}^T \mathbf{A}_{k_1 k_2}^{-1} \alpha_{k_2},$$

and we also used

$$\operatorname{tr} \mathbf{A}_{k_1 k_2}^{-1} \xi_{k_1 k_2} \leq \sigma^{-2} \operatorname{tr} \xi_{k_1 k_2} \leq \sigma^{-4}.$$

Repeating a step similar to (4.25), by Theorem 1 of [3] one can then conclude that

$$\begin{aligned}
&\frac{1}{N^4} \sum_{k_1 \neq k_2}^{K} E(\operatorname{tr} \mathbf{A}_{k_1 k_2}^{-2} - \operatorname{tr} E\mathbf{A}_{k_1 k_2}^{-2})^2 \\
&\leq \frac{MK^2}{N^4} + \frac{MK^2}{N^4} E(\operatorname{tr} \mathbf{A}^{-2} - \operatorname{tr} E\mathbf{A}^{-2})^2 \\
&\to 0,
\end{aligned} \quad (4.26)$$

as $N \to \infty$. Again, by an argument analogous to (4.25), one can find

$$\begin{aligned}
&E(\operatorname{tr} \mathbf{A}_{k_1 k_2}^{-1} - Nb_N)^4 \\
&\leq ME(\operatorname{tr} \mathbf{A}^{-1} - E \operatorname{tr} \mathbf{A}^{-1})^4 + M(E \operatorname{tr} \mathbf{A}^{-1} - Nb_N)^4 + M.
\end{aligned}$$

The second term on the right-hand side of the above inequality is bounded by the argument of Theorem 1 of [3]. We also claim that the first item on the right-hand side has an order $O(N)$. To see it, set $\mathcal{F}_j = \sigma(s_1, \ldots, s_j)$ and $E_j(\cdot) = E(\cdot|\mathcal{F}_j)$. By decomposing as the sum of a martingale difference sequence and using the Burkholder inequality, we have

$$\begin{aligned}
&E(\operatorname{tr} \mathbf{A}^{-1} - E \operatorname{tr} \mathbf{A}^{-1})^4 \\
&= E\left(\sum_{k=1}^{K} (E_k - E_{k-1}) \operatorname{tr}(\mathbf{A}^{-1} - \mathbf{A}_k^{-1})\right)^4
\end{aligned} \quad (4.27)$$

$$\leq ME\left(\sum_{k=1}^{K}((E_k - E_{k-1})\eta_k)^2\right)^2$$

$$\leq KM \sum_{k=1}^{K} E\left((E_k - E_{k-1})\frac{p_k \mathbf{s}_k^T \mathbf{A}_k^{-2} \mathbf{s}_k}{1 + p_k(1/N)\operatorname{tr}\mathbf{A}_k^{-1}}\right)^4$$

$$+ KM \sum_{k=1}^{K} E\left((E_k - E_{k-1})\frac{\eta_k(p_k \hat{s}_k - (p_k/N)\operatorname{tr}\mathbf{A}_k^{-1})}{1 + (p_k/N)\operatorname{tr}\mathbf{A}_k^{-1}}\right)^4$$

$$\leq KM \sum_{k=1}^{K} E\left(\mathbf{s}_k^T E_k\left(\frac{\mathbf{A}_k^{-2}}{1 + (p_k/N)\operatorname{tr}\mathbf{A}_k^{-1}}\right)\mathbf{s}_k\right.$$

$$\left. - \frac{1}{N} E_k \frac{\operatorname{tr}\mathbf{A}_k^{-2}}{1 + (p_k/N)\operatorname{tr}\mathbf{A}_k^{-1}}\right)^4$$

$$+ \frac{KM}{\sigma^8} \sum_{k=1}^{K} E\left(\hat{s}_k - \frac{1}{N}\operatorname{tr}\mathbf{A}_k^{-1}\right)^4$$

$$\leq \frac{MK^2}{N} = O(N),$$

where

$$\eta_k = \frac{p_k \mathbf{s}_k^T \mathbf{A}_k^{-2} \mathbf{s}_k}{1 + p_k \mathbf{s}_k^T \mathbf{A}_k^{-1} \mathbf{s}_k},$$

and we also use (4.18) and the equality

$$\eta_k = \frac{p_k \mathbf{s}_k^T \mathbf{A}_k^{-2} \mathbf{s}_k}{1 + (p_k/N)\operatorname{tr}\mathbf{A}_k^{-1}} - \frac{\eta_k(p_k \hat{s}_k - (p_k/N)\operatorname{tr}\mathbf{A}_k^{-1})}{1 + (p_k/N)\operatorname{tr}\mathbf{A}_k^{-1}}.$$

Combining the above one can conclude that

$$\frac{1}{N^4} \sum_{k_1 \neq k_2}^{K} E(\operatorname{tr}\mathbf{A}_{k_1 k_2}^{-1} - Nb_N)^4 \leq \frac{M}{N} \to 0 \quad \text{as } N \to \infty.$$

The basic inequality $(E|X|)^2 \leq EX^2$ implies the remaining term in (4.24) also goes to 0 as $N \to \infty$. Hence $U_{121}$ can be ignored.

Consider the term $U_{122}$ second. By (4.23) we have

$$(4.28) \quad \begin{aligned} & E((p_{k_2}\alpha_{k_2}\mathbf{s}_{k_1}^T \mathbf{A}_{k_1 k_2}^{-1}\mathbf{s}_{k_2}\mathbf{s}_{k_2}^T \mathbf{A}_{k_1 k_2}^{-1}\mathbf{s}_{k_1})^2 | \mathbf{s}_{k_2}, \mathbf{A}_{k_1 k_2}^{-1}) \\ & = \frac{3}{N^2}(p_{k_2}\alpha_{k_2}\mathbf{s}_{k_2}^T \mathbf{A}_{k_1 k_2}^{-2}\mathbf{s}_{k_2})^2. \end{aligned}$$

It follows that

$$U_{122} = \sum_{k_1 \neq k_2}^{K} E[E(p_{k_1 k_2} \zeta_{k_1 k_2} | \mathbf{s}_{k_2}, \mathbf{A}_{k_1 k_2}^{-1}) \beta_{k_2 k_1}]$$

$$\leq \frac{M}{N^2} \sum_{k_1 \neq k_2}^{K} E(\mathbf{s}_{k_2}^T \mathbf{A}_{k_1 k_2}^{-1} \mathbf{s}_{k_2} - b_N)^2$$

(4.29)
$$\leq \frac{MK^2}{N^3} + \frac{M}{N^2} \sum_{k_1 \neq k_2}^{K} E\left(\frac{1}{N} \operatorname{tr} \mathbf{A}_{k_1 k_2}^{-1} - b_N\right)^2$$

$$\leq \frac{MK^2}{N^3} + \frac{MK^2}{N^4} E(\operatorname{tr} \mathbf{A}^{-1} - N b_N)^2$$

$$= O\left(\frac{1}{N}\right).$$

Similarly, the term $U_{123}$ converges to zero.

Third, consider the term $U_{124}$. To simplify the notation, we write

$$s_{k_1}^{(j)} = \mathbf{s}_{k_1}^T \mathbf{A}_{k_1 k_2}^{-j} \mathbf{s}_{k_1}, \qquad s_{k_2}^{(j)} = \mathbf{s}_{k_2}^T \mathbf{A}_{k_1 k_2}^{-j} \mathbf{s}_{k_2}, \qquad j = 1, 2,$$

$$\hat{s}_{k_1 k_2} = \mathbf{s}_{k_1}^T \mathbf{A}_{k_1 k_2}^{-1} \mathbf{s}_{k_2} \mathbf{s}_{k_2}^T \mathbf{A}_{k_1 k_2}^{-1} \mathbf{s}_{k_1}.$$

According to (4.28) one can find that

$$U_{124} = \sum_{k_1 \neq k_2}^{K} E(p_{k_1 k_2} \alpha_{k_2}^2 \gamma_{k_1 k_2}^2 \gamma_{k_1 k_2}^2 \alpha_{k_1}^2)$$

$$- \frac{9}{N^4} \sum_{k_1 \neq k_2}^{K} p_{k_1}^2 p_{k_2}^2 p_{k_1 k_2} E(p_{k_2} \alpha_{k_2} s_{k_2}^{(2)})^2 E(p_{k_1} \alpha_{k_1} s_{k_1}^{(2)})^2$$

(4.30)
$$\leq M \sum_{k_1 \neq k_2}^{K} E\alpha_{k_2}^2 (p_{k_2} \hat{s}_{k_1 k_2})^4 + \frac{MK^2}{N^4}$$

$$\leq M \varepsilon_N = o(1).$$

Fourth, since the composition of the terms $U_{125}$ and $U_{127}$ is similar, we analyze only the $U_{125}$ term. From (4.28) we obtain

(4.31) $$U_{125} = -2 \sum_{k_1 \neq k_2}^{K} E(p_{k_1 k_2} \alpha_{k_2}^2 \alpha_{k_1} \gamma_{k_1 k_2}^2 p_{k_2} \gamma_{k_1 k_2} (s_{k_2}^{(1)} - b_N))$$

(4.32) $$+ \frac{3}{N^2} \sum_{k_1 \neq k_2}^{K} p_{k_1 k_2} E(p_{k_1} p_{k_2} s_{k_2}^{(2)} \alpha_{k_2})^2 E(p_{k_2} \gamma_{k_1 k_2} (s_{k_2}^{(1)} - b_N) \alpha_{k_1}).$$

By the Hölder inequality and (4.28) one can find

$$|E(p_{k_2}\gamma_{k_1k_2}(s_{k_2}^{(1)} - b_N)\alpha_{k_1})|$$
$$(4.33) \qquad \leq M(E(p_{k_1}\hat{s}_{k_1k_2}\alpha_{k_1})^2)^{1/2}(E(s_{k_2}^{(1)} - b_N)^2)^{1/2}$$
$$\leq \frac{M}{N^{3/2}},$$

which implies that the term in (4.31) converges to zero. Moreover, the absolute value of each summand in (4.31) is not larger than

$$ME|\alpha_{k_2}^2(p_{k_2}\hat{s}_{k_1k_2})^3(s_{k_2}^{(1)} - b_N)|$$
$$= ME[|s_{k_2}^{(1)} - b_N|E(|p_{k_2}\hat{s}_{k_1k_2}|^3\alpha_{k_2}^2|\mathbf{A}_{k_1k_2}^{-1},\mathbf{s}_{k_2})]$$
$$\leq ME\left[|s_{k_2}^{(1)} - b_N|E\left(\left|p_{k_2}\hat{s}_{k_1k_2} - \frac{1}{N}p_{k_2}s_{k_2}^{(2)}\right|^3\alpha_{k_2}^2\Big|\mathbf{A}_{k_1k_2}^{-1},\mathbf{s}_{k_2}\right)\right]$$
$$+ \frac{M}{N^3}E[|s_{k_2}^{(1)} - b_N|E(|p_{k_2}s_{k_2}^{(2)}|^3\alpha_{k_2}^2|\mathbf{A}_{k_1k_2}^{-1},\mathbf{s}_{k_2})]$$
$$\leq \frac{M(E|v_{11}|^6 + 2)}{N^3}E|(s_{k_2}^{(1)} - b_N)s_{k_2}^{(2)}|$$
$$\leq \frac{M}{N^{5/2}},$$

which leads the sum in (4.31) to converge to zero. So $U_{125}$ converges to zero as $N \to \infty$.

Fifth, consider $U_{126}$ ($U_{128}$ can be analyzed similarly):

$$U_{126}$$
$$(4.34) \qquad = -2\sum_{k_1 \neq k_2}^{K} p_{k_1k_2}E[p_{k_1}^2(s_{k_1}^{(1)} - b_N)^2 p_{k_2}\gamma_{k_1k_2}(s_{k_2}^{(1)} - b_N)\alpha_{k_1}]$$
$$(4.35) \qquad + 2\sum_{k_1 \neq k_2}^{K} p_{k_1k_2}p_{k_1}^2 E(s_{k_1}^{(1)} - b_N)^2 Ep_{k_2}\gamma_{k_1k_2}(s_{k_2}^{(1)} - b_N)\alpha_{k_1}.$$

From the estimate (4.33), the sum in (4.35) has an order $O(N^{-1/2})$. On the other hand, each summand in the sum in (4.34) can be rewritten as

$$E[p_{k_1}^2 p_{k_2} p_{k_1k_2}(s_{k_1}^{(1)} - b_N)^2 \alpha_{k_1} E(\gamma_{k_1k_2}(s_{k_2}^{(1)} - b_N)|\mathbf{A}_{k_1k_2}^{-1},\mathbf{s}_{k_1})]$$
$$\leq E[p_{k_1}^2 \alpha_{k_1}(s_{k_1}^{(1)} - b_N)^2 (E(\gamma_{k_1k_2}^2|\mathbf{A}_{k_1k_2}^{-1},\mathbf{s}_{k_1}))^{1/2}$$
$$\times (E((s_{k_2}^{(1)} - b_N)^2|\mathbf{A}_{k_1k_2}^{-1},\mathbf{s}_{k_1}))^{1/2}]$$

$$\leq \frac{M}{N} E[p_{k_1}^2 (s_{k_1}^{(1)} - b_N)^2 p_{k_1} s_{k_1}^{(2)} \alpha_{k_1} (E((s_{k_2}^{(1)} - b_N)^2 | \mathbf{A}_{k_1 k_2}^{-1}, \mathbf{s}_{k_1}))^{1/2}]$$

$$\leq \frac{M}{N^{5/2}},$$

and, thus, the sum in (4.34) converges to zero. Thus $U_{126}$ converges to zero as well.

Finally, consider the term $U_{129}$. As was done for the terms $U_{124}$, $U_{125}$ and $U_{126}$, the $U_{129}$ term is split into the sum of two terms. It follows from (4.33) that one of them,

$$4 \sum_{k_1 \neq k_2}^{K} p_{k_1 k_2} E(p_{k_1} \gamma_{k_1 k_2} (s_{k_1}^{(1)} - b_N) \alpha_{k_2}) E(p_{k_2} \gamma_{k_1 k_2} (s_{k_2}^{(1)} - b_N) \alpha_{k_1}),$$

converges to zero. The other term is

$$4 \sum_{k_1 \neq k_2}^{K} E(p_{k_1 k_2} p_{k_1} \gamma_{k_1 k_2} (s_{k_1}^{(1)} - b_N) \alpha_{k_2} p_{k_2} \gamma_{k_1 k_2} (s_{k_2}^{(1)} - b_N) \alpha_{k_1}).$$

The absolute value of each of the above summands is not larger than

$$(E(p_{k_1} \hat{s}_{k_1 k_2} (s_{k_1}^{(1)} - b_N) \alpha_{k_1})^2)^{1/2} \times (E(p_{k_2} \hat{s}_{k_1 k_2} (s_{k_2}^{(1)} - b_N) \alpha_{k_2})^2)^{1/2}.$$

Note that

$$E(p_{k_1} \hat{s}_{k_1 k_2} (s_{k_1}^{(1)} - b_N) \alpha_{k_1})^2$$
$$= E(((s_{k_1}^{(1)} - b_N) \alpha_{k_1})^2 E((p_{k_1} \hat{s}_{k_1 k_2})^2 | \mathbf{A}_{k_1 k_2}^{-1}, \mathbf{s}_{k_1}))$$
$$= \frac{3}{N^2} E((s_{k_1}^{(1)} - b_N) p_{k_1} s_{k_1}^{(2)} \alpha_{k_1})^2$$
$$\leq \frac{3}{\sigma^4 N^2} E(s_{k_1}^{(1)} - b_N)^2 \leq \frac{M}{N^3}.$$

Similarly,

$$E(p_{k_2} \hat{s}_{k_1 k_2} (s_{k_2}^{(1)} - b_N) \alpha_{k_2})^2 \leq M/N^3.$$

Hence $U_{129}$ converges to zero.

Summarizing the above argument, one can conclude that the variance of the term $U_1$ converges to zero as $N \to \infty$ and thus, it is sufficient to compute the asymptotic value of its expectation, which can be accomplished as follows:

$$\sum_{k=1}^{K} \hat{p}_k E(\hat{s}_k - b_N)^2$$

$$(4.36) \quad = \sum_{k=1}^{K} \hat{p}_k E\left(\hat{s}_k - \frac{1}{N} \operatorname{tr} \mathbf{A}_k^{-1}\right)^2 + \sum_{k=1}^{K} \hat{p}_k E\left(\frac{1}{N} \operatorname{tr} \mathbf{A}_k^{-1} - b_N\right)^2$$

$$= \frac{2}{N^2} \sum_{k=1}^{K} \hat{p}_k E \operatorname{tr} \mathbf{A}_k^{-2} + o(1)$$

$$= \frac{2}{N^2} E \operatorname{tr} \mathbf{A}^{-2} \sum_{k=1}^{K} \hat{p}_k + o(1)$$

$$\to 2c \int \frac{dF^{c,H}(x)}{(x+\sigma^2)^2} \int \frac{x^2}{(1+xb)^3} dH(x),$$

where $\hat{p}_k = p_k^2/(1+b_N p_k)^3$, and in the second and third equalities we use a trick similar to (4.25).

Note

$$(4.37) \quad F^{c,H} = (1-c)I_{[0,\infty)} + c\underline{F}^{c,H},$$

where $\underline{F}^{c,H}$ represents the limiting spectral distribution of $\mathbf{P}^{1/2}\mathbf{S}^T\mathbf{S}\mathbf{P}^{1/2}$ with $\mathbf{P}^{1/2} = \operatorname{diag}(\sqrt{p_1}, \ldots, \sqrt{p_K})$. From (4.37), one can get

$$\operatorname{tr} \mathbf{A}^{-1} - Nb_N = \operatorname{tr}(\mathbf{P}^{1/2}\mathbf{S}^T\mathbf{S}\mathbf{P}^{1/2} + \sigma^2 I)^{-1} - K \int \frac{d\underline{F}^{c_N, H_N}(x)}{x+\sigma^2}$$

and Theorem 1.1 of [3] is then applicable. Thus, for $U_2$ one has a central limit theorem and it then suffices to show that $U_5$ converges to zero in probability. It is obvious that the term $U_5$ becomes zero when $p_1 = \cdots = p_K$; however, its convergence in probability appears to be somewhat troublesome when the powers of the users are not the same. We will provide an abridged analysis for this case. Set

$$a_k = \frac{1}{(1+p_k b_N)^2}, \qquad a = \int \frac{dH_N(x)}{(1+xb_N)^2}.$$

Using steps analogous to (4.25), one can obtain

$$U_5 = \sum_{k=1}^{K} p_k \left(\hat{s}_k - \frac{1}{N} \operatorname{tr} \mathbf{A}_k^{-1}\right)(a_k - a)$$

$$+ (\operatorname{tr} \mathbf{A}^{-1} - Nb_N)\frac{1}{N} \sum_{k=1}^{K} p_k(a_k - a) + o_p(1)$$

$$= \sum_{k=1}^{K} p_k \left(\hat{s}_k - \frac{1}{N} \operatorname{tr} \mathbf{A}_k^{-1}\right)(a_k - a) + o_p(1).$$

Hence it is necessary to show that

$$(4.38) \quad \hat{U}_5 \triangleq E\left(\sum_{k=1}^{K} p_k\left(\hat{s}_k - \frac{1}{N}\operatorname{tr}\mathbf{A}_k^{-1}\right)(a_k - a)\right)^2$$

converges to zero. Expanding out the right-hand side of (4.38) one can get

$$\hat{U}_5 = U_{51} + U_{52},$$

where

$$U_{51} = \sum_{k=1}^{K} E\left(p_k\left(\hat{s}_k - \frac{1}{N}\operatorname{tr}\mathbf{A}_k^{-1}\right)(a_k - a)\right)^2,$$

$$U_{52} = \sum_{k_1 \neq k_2}^{K} E\bigg(p_{k_1}p_{k_2}\left(s_{k_1}^{(k_2)} - \frac{1}{N}\operatorname{tr}\mathbf{A}_{k_1}^{-1}\right)$$

$$\times \left(s_{k_2}^{(k_1)} - \frac{1}{N}\operatorname{tr}\mathbf{A}_{k_2}^{-1}\right)(a_{k_1} - a)(a_{k_2} - a)\bigg).$$

It is easy to see that

$$U_{51} \leq \frac{M}{N} \sum_{k=1}^{K} (p_k(a_k - a))^2 \to 0.$$

Regarding the term $U_{52}$, one can show that it converges to zero by an argument similar to that used for the preceding term $U_{12}$ and since the process is somewhat tedious, it is omitted.

For the computation of (1.11) and (1.12), without loss of generality, suppose $p = 1$; otherwise replace $\sigma^2$ by $\sigma^2/p$. As for the formulas (1.9) and (1.11) one can refer to, respectively, (1.18) and (5.13) of [3].

Now let us derive (1.12). It is shown in [7] that $m(z) = m_{F^{c,H}}(z)$, for each $z \in \mathbb{C}^+$, is the unique solution in $\mathbb{C}^+$ to the equation

$$(4.39) \quad m = -\left(z - c\int \frac{t\,dH(t)}{1+tm}\right)^{-1}.$$

From this equation, the inverse function has an explicit form

$$(4.40) \quad z = -\frac{1}{m} + c\int \frac{t\,dH(t)}{1+tm}$$

and one can then find for $H(t) = I_{[1,\infty)}(t)$

$$(4.41) \quad z = -\frac{1}{m(z)} + \frac{c}{1+m(z)}.$$

Suppose the $m_2$ contour encloses the $m_1$ contour (see [8] or [3] for the range of $m(x)$ for a real $x$ and contour of $m$). For a fixed $m_2$ it follows from (4.41) and the Cauchy residue theorem that

$$\int \frac{dm_1}{(z(m_1)+\sigma^2)(m_1-m_2)^2} = \frac{1}{\sigma^2}\int \frac{(m_1^2+m_1)\,dm_1}{(m_1-m_a)(m_1-m_b)(m_1-m_2)^2}$$

$$= \frac{2\pi i}{\sigma^2}\frac{m_b^2+m_b}{(m_b-m_2)^2(m_b-m_a)},$$

where

$$m_a = \frac{-(1+(c-1)/\sigma^2)+\sqrt{(1+(c-1)/\sigma^2)^2+4/\sigma^2}}{2},$$

$$m_b = \frac{-(1+(c-1)/\sigma^2)-\sqrt{(1+(c-1)/\sigma^2)^2+4/\sigma^2}}{2}.$$

Consequently,

$$\rho = \frac{1}{\sigma^2\pi i}\int \frac{m_b^2+m_b}{(m_b-m_2)^2(m_b-m_a)(z(m_2)+\sigma^2)}\,dm_2$$

$$= \frac{m_b^2+m_b}{\sigma^4\pi(m_b-m_a)i}\int \frac{(m_2^2+m_2)\,dm_2}{(m_2-m_b)^3(m_2-m_a)}$$

$$= \frac{2(m_b^2+m_b)(m_a^2+m_a)}{\sigma^4(m_b-m_a)^4} = \frac{2c}{((\sigma^2+c-1)^2+4\sigma^2)^2}.$$

**5. Proof of Corollary 1.1.** Using the Taylor expansion, one can find

$$\sum_{k=1}^{K}(\log(1+\beta_k)-\log(1+b_Np_k))$$

$$= \sum_{k=1}^{K}\frac{\beta_k-b_Np_k}{1+b_Np_k} - \sum_{k=1}^{K}\frac{(\beta_k-b_Np_k)^2}{2(1+b_Np_k)^2} + \sum_{k=1}^{K}\frac{(\beta_k-b_Np_k)^3}{3(1+\xi_k)^3}$$

$$= \int \frac{dH_N(x)}{1+xb_N}\sum_{k=1}^{K}(\beta_k-b_Np_k) - \sum_{k=1}^{K}\frac{(\beta_k-b_Np_k)^2}{2(1+b_Np_k)^2} + \sum_{k=1}^{K}\frac{(\beta_k-b_Np_k)^3}{3(1+\xi_k)^3}$$

$$+ \sum_{k=1}^{K}(\beta_k-b_Np_k)\left(\frac{1}{1+b_Np_k}-\int\frac{dH_N(x)}{1+xb_N}\right),$$

where each $\xi_k$ is located in the interval $[\beta_k,b_Np_k]$. Since

$$\sum_{k=1}^{K}\frac{(\beta_k-b_Np_k)^3}{3(1+\xi_k)^3} \leq \sum_{k=1}^{K}|\beta_k-b_Np_k|^3,$$

Corollary 1.1 holds by the argument of Theorem 1.3.

**Acknowledgments.** The authors would like to thank the Associate Editor and the referees for their helpful suggestions.

G.-M. Pan
M.-H. Guo
Department of Applied Mathematics
National Sun Yat Sen University
Taiwan
E-mail: stapgm@gmail.com
  guomh@math.nsysu.edu.tw

W. Zhou
Department of Statistics
and Applied Probability
National University of Singapore
Singapore 117546
E-mail: stazw@nus.edu.sg